\documentclass{elsart}
\usepackage{graphics}
\usepackage{graphicx}
\usepackage{epsfig}
\usepackage{amssymb}

\begin{document}

\begin{frontmatter}
\title{A Graph-Theoretic Approach to a Partial Order of Knots and Links}
\author{Toshiki Endo}
\ead{end@prf.jiyu.ac.jp}
\address{Research Institute for Mathematics and Information Science, Jiyu Gakuen College, Gakuencho 1-8-15, Higashikurume-shi, Tokyo 203-8521, Japan}
\author{Tomoko Itoh}
\ead{tomokoitoh@prf.jiyu.ac.jp}
\address{Research Institute for Mathematics and Information Science, Jiyu Gakuen College, Gakuencho 1-8-15, Higashikurume-shi, Tokyo 203-8521, Japan}
\author{Kouki Taniyama}
\ead{taniyama@waseda.jp}
\address{Department of Mathematics, School of Education, Waseda University, Nishi-Waseda 1-6-1, Shinjuku-ku, Tokyo 169-8050, Japan}

\maketitle

\begin{abstract}
We say that a link $L_1$ is an s-major of a link $L_2$ if any diagram of $L_1$ can be transformed into a diagram of $L_2$ by changing some crossings and smoothing some crossings. This relation is a partial ordering on the set of all prime alternating links. We determine this partial order for all prime alternating knots and links with the crossing number less than or equal to six. The proofs are given by graph-theoretic methods.
\end{abstract}


\end{frontmatter}

\section{Introduction}

Throughout this paper we work in the piecewise linear category. We use standard terminology and notation of knot theory, see for example \cite{A} and \cite{R}, and graph theory, see for example \cite{B} and \cite{D}. In particular we denote the complete graph on $n$ vertices by $K_n$, the $n$-cycle by $C_n$ and a graph on two vertices and $n$ multiple edges joining them by $\theta_n$.

We assume that all links are unordered, unoriented and contained in the 3-sphere ${\mathbb S}^3$. By a {\it link diagram}, or simply a {\it diagram}, we mean a regular diagram of a link in the 2-sphere ${\mathbb S}^2 \subset {\mathbb S}^3$. Note that a diagram has only finitely many transversal double points each of which has over/under crossing information. We call such a double point a {\it crossing}. A diagram without over/under crossing information is called a {\it projection}.

We say that two links $L_1$ and $L_2$ are {\it equivalent}, denoted by $L_1=L_2$, if there exists a possibly orientation reversing homeomorphism of ${\mathbb S}^3$ onto itself which maps $L_1$ to $L_2$. The equivalence class is called a link type. We do not distinguish between a link and its link type so long as no confusion occurs.

\begin{figure}[http]
\begin{center}
\includegraphics[width=140mm]{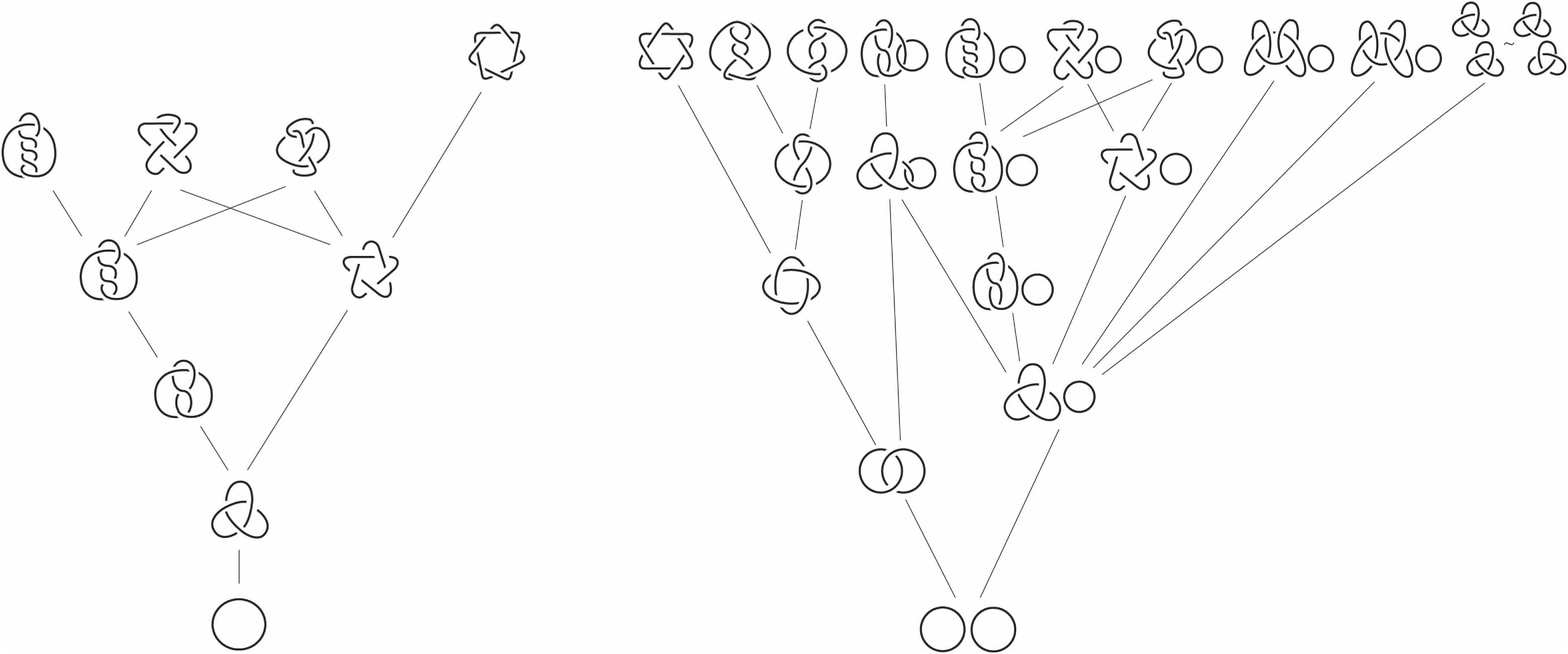}
\put(-261,150){$7_1$}
\put(-380,125){$6_1$} \put(-344,125){$6_2$} \put(-310,125){$6_3$}
\put(-360,92){$5_2$} \put(-290,92){$5_1$}
\put(-340,62){$4_1$}
\put(-325,30){$3_1$}
\put(-325,3){$0_1$}
\caption{Pre-orders of knots and 2-component links}
\label{fig0101}
\end{center}
\end{figure}

In \cite{K1} the third author defined a pre-ordering on the set of $\mu$-component links as follows. Let $L_1$ and $L_2$ be $\mu$-component links. Then we say that $L_1$ is a {\it major} of $L_2$, denoted by $L_1\geq L_2$, if every projection of $L_1$ is also a projection of $L_2$. In other words every diagram of $L_1$ can be transformed into a diagram of $L_2$ by changing over/under information at some crossings of the diagram of $L_1$. Then we also say that $L_2$ is a {\it minor} of $L_1$. The third author studied this order for knots in \cite{K1} and for 2-component links in \cite{K2}, and obtained two Hasse diagrams shown in Fig. \ref{fig0101}, where each line segment means that the upper one is a major of the lower one. In the sequel, the number representing a link is due to the Rolfsen's knot table in \cite{A}. 

We remark here that after \cite{K1} and \cite{K2} some related works are done. They are for example \cite{K-M}, \cite{T-Y}, \cite{K3}, \cite{H-J-O1}, \cite{H-J-O2}, \cite{H-T} and \cite{Nikkuni}. We also note that an application of the results in \cite{K1} and \cite{K2} to link signature are done in the forthcoming paper \cite{P-T}.

\begin{figure}[http]
\begin{center}
\includegraphics[width=40mm]{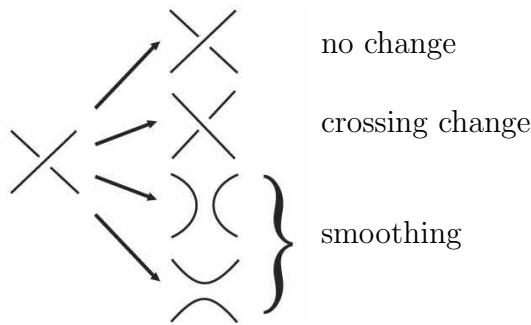}
\put(5,102){no change}
\put(5,72){crossing change}
\put(5,30){smoothing}
\caption{Four operations at a crossing point}
\label{fig0102}
\end{center}
\end{figure}

Note that this order is defined only to links with the same number of components, since any crossing change does not increase or decrease the number of link components. We now define an extended version of this order, which enables us to compare links with different numbers of components. The allowable operation connecting them is a {\it smoothing} operation at a crossing point, shown in Fig. \ref{fig0102}. 

\begin{figure}[http]
\begin{center}
\includegraphics[width=40mm]{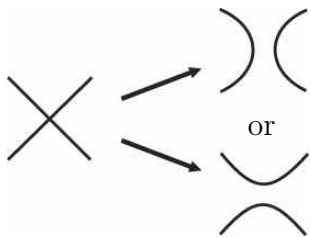}
\put(-23,38){or}
\caption{Smoothing a crossing of a projection}
\label{fig0102-2}
\end{center}
\end{figure}

{\bf Definition 1.1.}\quad
Let $L_1$ and $L_2$ be links. We say that $L_1$ is an {\it s-major} of $L_2$ if every diagram of $L_1$ can be transformed into a diagram of $L_2$ by applying one of the four operations illustrated in Fig. \ref{fig0102} at each crossing point of the diagram of $L_1$. We denote it by $L_1 \succeq L_2$. Then we also say that $L_2$ is an {\it s-minor} of $L_1$ and denote it by $L_2 \preceq L_1$. We call this order {\it smoothing order}. We note that $L_1$ is an s-major of $L_2$ if and only if every projection of $L_1$ can be transformed into a projection of $L_2$ by smoothing some crossing points of the projection of $L_1$ as illustrated in Fig. \ref{fig0102-2}.

By definition we have the following proposition.

{\bf Proposition 1.2.}\quad
{\it Let $L_1$ and $L_2$ be links of the same number of components. If $L_1\geq L_2$ then $L_1 \succeq L_2$.}

{\bf Historical remark.}\quad
This order has first defined by the third author and some results are announced without proof in \cite{K3}. In 2007 the second author re-defined this order without knowledge of \cite{K3}. We note that all results announced in \cite{K3} are contained in this paper.

Let ${\cal L}$ be the set of all links. The following two propositions immediately follow from the definition.

{\bf Proposition 1.3.}\quad
{\it The pair $({\cal L},\succeq)$ is a pre-ordered set. Namely, for any $L_1$, $L_2$ and $L_3$ in ${\cal L}$ the following (1) and (2) hold.

(1)  $L_1\succeq L_1$  (reflexive law).

(2) If $L_1 \succeq L_2$ and $L_2 \succeq L_3$, then $L_1 \succeq L_3$ (transitive law).}

{\bf Proposition 1.4.}\quad
{\it Let $L_1$ and $L_2$ be links. Suppose that $L_1\succeq L_2$. Then we have $c(L_1)\geq c(L_2)$ where $c(L)$ denotes the minimal crossing number of $L$.}

A link is said to be {\it prime} if it is non-splittable and every 2-sphere in ${\mathbb S}^3$ meeting the link transversely in two points bounds a trivial ball-arc pair. Thus we treat a trivial knot as a prime knot in this paper. It is known that every link is decomposed into finite number of prime links \cite{H}. These prime links are called {\it prime factors} of the link. A link is called alternating if it has a diagram in which over-crossing and under-crossing appear alternately. We denote the set of all prime alternating links by $\cal{PAL}$. Then we have the following proposition.

{\bf Proposition 1.5.}\quad
{\it The pair $(\cal{PAL}, \succeq)$ is a partially ordered set. Namely, in addition to the reflective law and the transitive law the following holds.

(3) If $L_1,L_2\in\cal{PAL}$, $L_1 \succeq L_2$ and $L_2 \succeq L_1$, then $L_1=L_2$ (the antisymmetric law).}

{\it Proof.}\quad
Suppose that $L_1,L_2\in\cal{PAL}$, $L_1 \succeq L_2$ and $L_2 \succeq L_1$. Then by Proposition 1.3 we have $c(L_1)=c(L_2)$. Let $\tilde L_1$ be a minimal crossing diagram of $L_1$. It is known that a minimal crossing diagram of a prime alternating link is always reduced alternating \cite{K} \cite{M} \cite{T}. Since smoothing decreases the number of crossings we have that a diagram $\tilde L_2$ of $L_2$ is obtained by changing some crossings of $\tilde L_1$. Then $\tilde L_2$  is also a minimal diagram of $L_2$. Therefore $\tilde L_2$ is also alternating. Since these diagrams are connected we have that they are either identical or differ by all crossings. The latter case implies that $L_2$ is a mirror image of $L_1$. Then by definition we have $L_1=L_2$. \hfill$\Box$

\begin{figure}[http]
\begin{center}
\includegraphics[width=60mm]{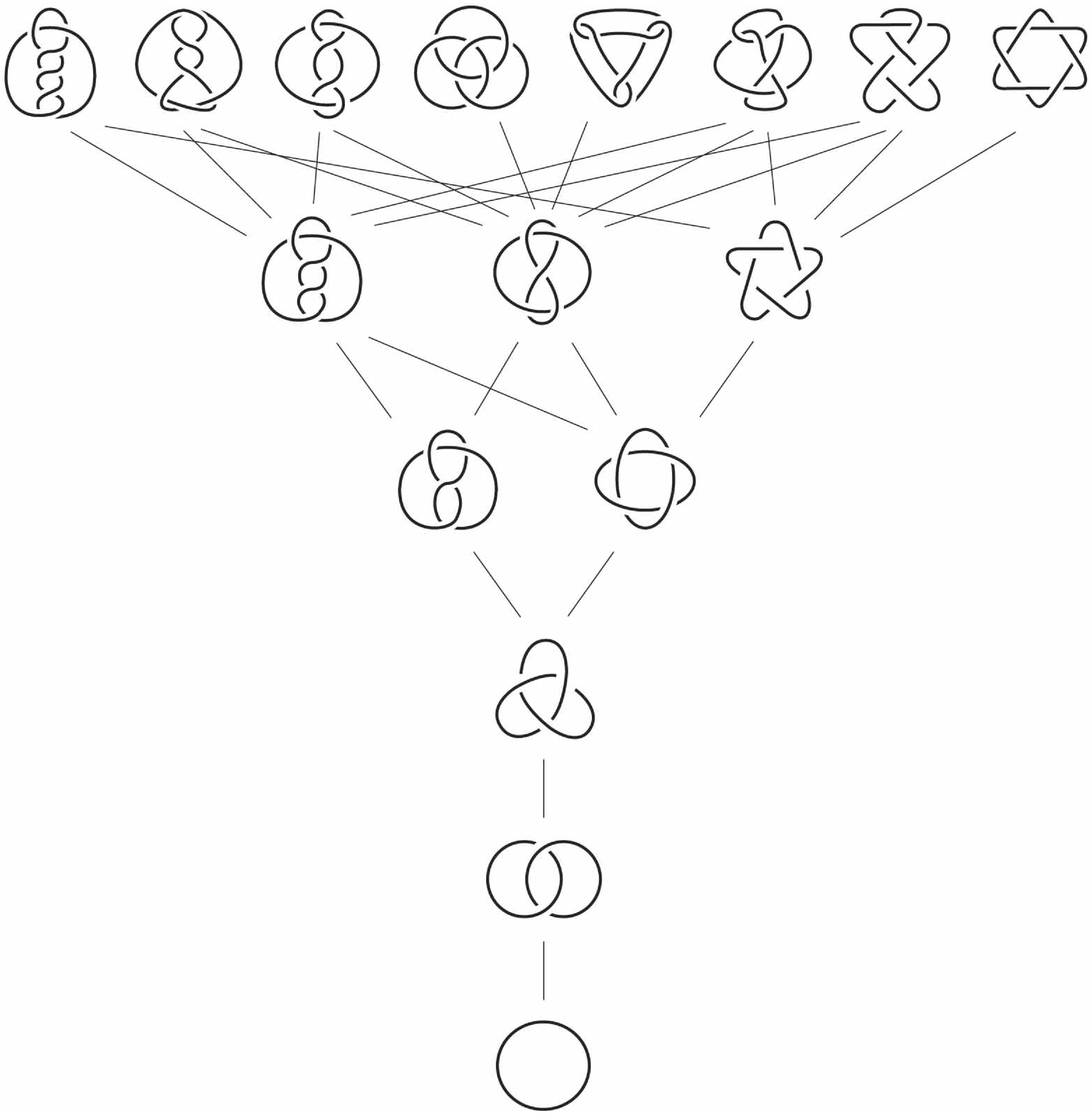}
\put(-166,178){$6_1$}\put(-144,178){$6^2_2$}\put(-122,178){$6^2_3$}\put(-100,178){$6^3_2$}\put(-78,178){$6^3_1$}\put(-56,178){$6_3$}\put(-34,178){$6_2$}\put(-12,178){$6^2_1$}
\put(-112,127){$5_2$} \put(-76,127){$5^2_1$} \put(-40,127){$5_1$}
\put(-91,93){$4_1$} \put(-61,93){$4^2_1$}
\put(-75,60){$3_1$}
\put(-75,33){$2^2_1$}
\put(-75,2){$0_1$}
\caption{Smoothing order of prime alternating knots and links with the crossing number $\leq 6$}
\label{fig0103}
\end{center}
\end{figure}

Our results are summarized by the Hasse diagram shown in Fig. \ref{fig0103}. Details are stated in the following section. 

Our strategy for proofs is graph-theoretic and is different from the methods in \cite{K1}\cite{K2}. It is well known that there is a correspondence between connected link diagrams on ${\mathbb S}^2$ and edge-signed plane graphs on ${\mathbb S}^2$. See for example \cite{A}. If we ignore the over/under crossing information we have a correspondence between connected link projections on ${\mathbb S}^2$ and plane graphs on ${\mathbb S}^2$. We briefly review this correspondence.

Let $L$ be a link, $\tilde L$ a connected diagram of $L$ on ${\mathbb S}^2$ and $\hat L$ its underlying projection. Then $\hat L$ is a connected 4-regular plane graph on ${\mathbb S}^2$. Let ${\cal C}$ be a coloring of the regions ${\mathbb S}^2-\hat L$ black and white such that adjacent regions have different colors. Note that there are two such colorings. Let $G(\hat L,{\cal C})$ be a plane graph on ${\mathbb S}^2$ contained in the closure of the union of black regions whose vertices lie in the black regions in one to one correspondence and whose edges are in one to one correspondence to the crossings of $\hat L$ so that each of them joins the the vertices in two black regions meeting at a crossing. 

Conversely, for a connected plane graph $G$ on ${\mathbb S}^2$, we take disks on ${\mathbb S}^2$ such that each disk contains just one vertex and two disks containing two adjacent vertices meets at the middle point of the edge joining them. Then the boundary of the union of  such disks is a 4-regular graph on ${\mathbb S}^2$ so that we may suppose it a link projection. We denote this link projection by $\hat L(G)$. Note that $\hat L(G(\hat L,{\cal C}))=\hat L$.

We say that a graph $H$ is a {\it minor} of a graph $G$ if $H$ is obtained from $G$ by a series of edge-contraction and taking subgraph. In addition we only consider edge-signed graphs on ${\mathbb S}^2$ and assume that edge contraction and taking subgraph preserves the signs of the survived edges and edge-contraction is performed on ${\mathbb S}^2$ so that it respects the embedding of the graph into ${\mathbb S}^2$.

The following proposition can be shown by standard arguments in graph theory. We omit the proof.

{\bf Proposition 1.6.}\quad
{\it Let $G$ and $H$ be connected graphs. Suppose that $H$ is a minor of $G$. Then there is a sequence of connected graphs $G=G_0,G_1,\cdots,G_n=H$ such that $G_{i+1}$ is obtained from $G_i$ by deleting an edge or by contracting an edge for each $i\in\{0,1,\cdots,n-1\}$.}

Then we immediately have the following proposition that is a key to prove the theorems of the next section.

{\bf Proposition 1.7.}\quad
{\it Let $\hat L$ be a connected link diagram. Let ${\cal C}$ be a coloring of ${\mathbb S}^2-\hat L$. Suppose that a connected graph $H$ is a minor of $G(\hat L,{\cal C})$. Then we have that the link projection $\hat L(H)$ can be obtained from $\hat L$ by smoothing some crossings of $\hat L$.}


\section{Determining the Smoothing Order of Prime Alternating Knots and Links up to Six Crossings}

\subsection{Links that are s-majors of the trivial knot}

{\bf Proposition 2.1.}\quad
{\it An s-major of a non-split link is non-splittable.}

{\it Proof.}\quad
A split link has a disconnected diagram from which no connected diagram arises. Therefore a split link cannot be an s-major of a non-splittable link. \hfill$\Box$

{\bf Theorem 2.2.}\quad
{\it A link is an s-major of a trivial knot if and only if it is non-splittable.} 

{\it Proof.}\quad
By Proposition 2.1 we have that a split link is not an s-major of a trivial knot. Let $L$ be a non-splittable link and $\hat L$ a projection of $L$. Then $\hat L$ is connected. Let ${\cal C}$ be a coloring of ${\mathbb S}^2-\hat L$. Then $K_1$ is a minor of $G(\hat L,{\cal C})$. Note that $\hat L(K_1)$ is a projection of a trivial knot.
Then by Proposition 1.7 we have that $L$ is a s-major of a trivial knot. \hfill$\Box$

\subsection{Links that are s-majors of the Hopf link}

\begin{figure}[http]
\begin{center}
\includegraphics[width=50mm]{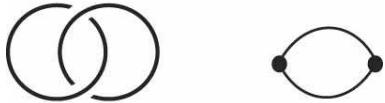}
\caption{The Hopf link and its graph}
\label{fig0202}
\end{center}
\end{figure}

A diagram of the Hopf link and a plane graph corresponding to its underlying projection are illustrated in Fig. \ref{fig0202}. 

We say that a crossing $c$ of a link projection $\hat L$ is {\it nugatory} if the number of the connected components of $\hat L-c$ is greater than that of $\hat L$.
We say that $\hat L$ is {\it reduced} if it has no nugatory crossings. We note that a nugatory crossing appears as a loop or a cut edge in the corresponding plane graph.
It is easy to see that we may only consider reduced projections. That is, we have the following proposition.

{\bf Proposition 2.3.}\quad
{\it Let $L_1$ and $L_2$ be links. Suppose that for any reduced projection $\hat L_1$ of $L_1$ there is a projection $\hat L_2$ of $L_2$ that is obtained from $\hat L_1$ by smoothing some crossings. Then $L_1$ is an s-major of $L_2$.} 

From now on all link projections are supposed to be reduced. 

{\bf Theorem 2.4.}\quad
{\it A link is an s-major of a Hopf link if and only if it is non-splittable and it is not a trivial knot.} 

{\it Proof.}\quad
The \lq only if' part follows from Proposition 2.1 and the fact that a trivial knot is not an s-major of a Hopf link. Let $L$ be a non-splittable link that is not a trivial knot and $\hat L$ a reduced projection of $L$. Let ${\cal C}$ be a coloring of ${\mathbb S}^2-\hat L$. Then $G(\hat L,{\cal C})$ is a connected graph without loops nor cut edges. Since $L$ is not a trivial knot we have that $G(\hat L,{\cal C})$ has at least two edges. Then we have that a 2-cycle is a minor of $G(\hat L,{\cal C})$. Then by Proposition 1.7 we have the result.\hfill$\Box$

\subsection{Links that are s-majors of the trefoil knot}

\begin{figure}[http]
\begin{center}
\includegraphics[width=70mm]{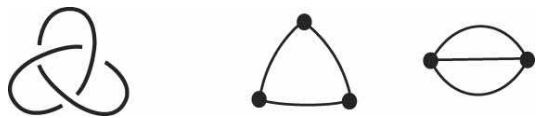}
\caption{The trefoil knot and its graphs}
\label{fig0203}
\end{center}
\end{figure}

A diagram of the trefoil knot and plane graphs $C_3$ and $\theta_3$ corresponding to its underlying projection are illustrated in Fig. \ref{fig0203}. 

{\bf Theorem 2.5.}\quad
{\it A link is an s-major of the trefoil knot if and only if it is non-splittable and it has a prime factor that is not a Hopf link.}

{\it Proof.}\quad
It is easily seen that a connected sum of some Hopf links is not an s-major of the trefoil knot. Therefore the \lq only if' part follows.
Let $L$ be a non-splittable link that has a prime factor which is not a Hopf link.
Let $\hat L$ be a reduced projection of $L$ and let ${\cal C}$ be a coloring of ${\mathbb S}^2-\hat L$. Then $G(\hat L,{\cal C})$ is a connected graph without loops nor cut edges. Note that $G(\hat L,{\cal C})$ may have some cut vertices. By the assumption there is a block $H$ of $G(\hat L,{\cal C})$ that corresponds to not necessarily one but some prime factors of $L$ such that at least one of them is not a Hopf link.
If $H$ has a cycle of length three or more then we have that the 3-cycle $C_3$ is a minor of $H$, hence of $G(\hat L,{\cal C})$. Suppose that the length of any cycle of $H$ is less than or equal to two. Then $H$ must be the graph on two vertices and three or more edges joining them. Then $\theta_3$ is a minor of $H$ and $G(\hat L,{\cal C})$. By Proposition 1.7 we have the result.\hfill$\Box$

\subsection{Links that are s-majors of the $(2,4)$-torus link}

\begin{figure}[http]
\begin{center}
\includegraphics[width=70mm]{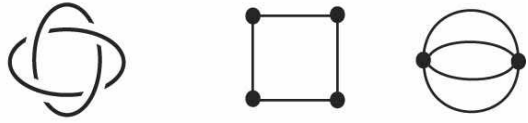}
\caption{The $(2,4)$-torus link and its graphs}
\label{fig0204}
\end{center}
\end{figure}

A diagram of the $(2,4)$-torus link and plane graphs $C_4$ and $\theta_4$ corresponding to its underlying projection are illustrated in Fig. \ref{fig0204}.

\begin{figure}[http]
\begin{center}
\includegraphics[width=45mm]{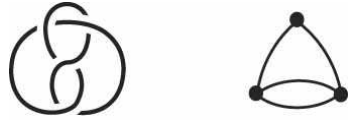}
\caption{The figure eight knot and its graph}
\label{fig02041}
\end{center}
\end{figure}

We also illustrate here in Fig. \ref{fig02041} a diagram of the figure eight knot $4_1$ and a plane graph corresponding to its underlying projection.

{\bf Theorem 2.6.}\quad
{\it A link is an s-major of the $(2,4)$-torus link if and only if it is non-splittable and it has a prime factor that is none of a Hopf link, a trefoil knot and a figure eight knot.}

{\it Proof.}\quad
It is easily seen that a connected sum of some Hopf links, trefoil knots and figure eight knots is not an s-major of the $(2,4)$-torus link. Therefore the \lq only if' part follows.
Let $L$ be a non-splittable link that has a prime factor which is none of a Hopf link, a trefoil knot and a figure eight knot.
Let $\hat L$ be a reduced projection of $L$ and let ${\cal C}$ be a coloring of ${\mathbb S}^2-\hat L$. By the assumption there is a block $H$ of $G(\hat L,{\cal C})$ that corresponds to not necessarily one but some prime factors of $L$ such that at least one of them is none of a Hopf link, a trefoil knot and a figure eight knot.
Then we have that $H$ has four or more edges. If $H$ has a cycle of length four or more, then we have that a 4-cycle $C_4$ is a minor of $H$ and hence of $G(\hat L,{\cal C})$. Thus, in this case, $L$ is an s-major of the $(2,4)$-torus link. Suppose that the length of any cycle of $H$ is less than or equal to three. If $H$ has two vertices connected by four or more internally disjoint paths, then we have that $\theta_4$ is a minor of $H$ and hence of $G(\hat L,{\cal C})$ and we have that $L$ is an s-major of the $(2,4)$-torus link.

\begin{figure}[http]
\begin{center}
\includegraphics[width=15mm]{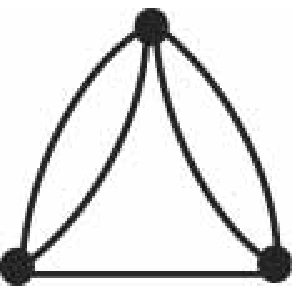}
\caption{}
\label{fig02041-2}
\end{center}
\end{figure}

Suppose that every pair of vertices of $H$ has at most three internally disjoint paths between them. Then we have that either $H$ contains the graph illustrated in Fig. \ref{fig02041-2} from which $\theta_4$ is obtained by an edge contraction, or $H$ is a proper minor of the graph in illustrated in Fig. \ref{fig02041-2}. Then it is easy to check that the graph only corresponds to a trivial knot, a trefoil knot or a figure eight knot. \hfill$\Box$

\subsection{Links that are s-majors of the figure eight knot}

{\bf Theorem 2.7.}\quad
{\it A link is an s-major of the figure eight knot if and only if it is non-splittable and it has a prime factor that is none of $(2,n)$-torus knots and links.}

{\it Proof.}\quad
It is easily seen that a connected sum of some $(2,n)$-torus knots and links is not an s-major of the figure eight knot. Therefore the \lq only if' part follows.
Let $L$ be a non-splittable link that has a prime factor which is none of $(2,n)$-torus knots and links.
Let $\hat L$ be a reduced projection of $L$ and let ${\cal C}$ be a coloring of ${\mathbb S}^2-\hat L$. By the assumption there is a block $H$ of $G(\hat L,{\cal C})$ that corresponds to not necessarily one but some prime factors of $L$ such that at least one of them is none of $(2,n)$-torus knots and links. Then we have that $H$ is neither a cycle nor a $\theta_n$-curve. Then it is easy to see that the graph illustrated in Fig. \ref{fig02041} is a minor of $H$. \hfill$\Box$

\begin{figure}[http]
\begin{center}
\includegraphics[width=40mm]{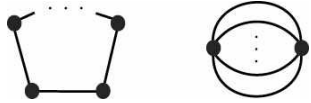}
\caption{Graphs corresponding to the $(2,n)$-torus link}
\label{fig0205}
\end{center}
\end{figure}

\subsection{Links that are s-majors of the Whitehead link}

\begin{figure}[http]
\begin{center}
\includegraphics[width=70mm]{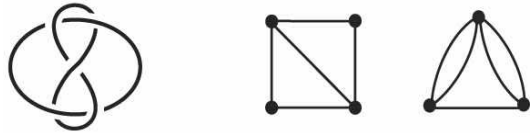}
\caption{The Whitehead link and its graphs}
\label{fig0206}
\end{center}
\end{figure}

A diagram of the Whitehead link $5^2_1$ and plane graphs corresponding to its underlying projection are illustrated in Fig. \ref{fig0206}.
One of them is a 4-cycle with one diagonal, and the other is the dual which is a 3-cycle with two multiple edges.

{\bf Theorem 2.8.}\quad
{\it A link is an s-major of the Whitehead link if and only if it is non-splittable and it has a prime factor that is none of $(2,n)$-torus knots and links and twist knots.}

{\it Proof.}\quad
It is easily seen that a connected sum of some $(2,n)$-torus knots and links and twist knots is not an s-major of the Whitehead link. Therefore the \lq only if' part follows.
Let $L$ be a non-splittable link that has a prime factor which is none of $(2,n)$-torus knots and links and twist knots. Let $\hat L$ be a reduced projection of $L$ and let ${\cal C}$ be a coloring of ${\mathbb S}^2-\hat L$. By the assumption there is a block $H$ of $G(\hat L,{\cal C})$ that corresponds not necessarily one but some prime factors of $L$ such that at least one of them is none of $(2,n)$-torus knots and links and twist knots.
If $H$ has a cycle of length greater than or equal to four and the cycle has a diagonal path, or if $H$ has a cycle of length three or more such that at least two edges of the cycle are multiple edges, then we obtain a graph illustrated in Fig. \ref{fig0206}. 

\begin{figure}[http]
\begin{center}
\includegraphics[width=46mm]{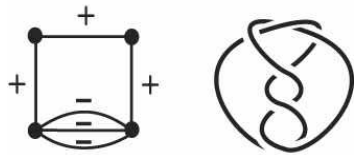}
\caption{A graph corresponding to a Whitehead link}
\label{fig02062}
\end{center}
\end{figure}

If $H$ has a cycle of length greater than or equal to four with multiple three edges, then we get the graph in Fig. \ref{fig02062}, which also corresponds to a Whitehead link. 

Thus, we may assume that if $H$ has a cycle of length greater than or equal to four then the cycle has no diagonals, no multiple three edges, and no distinct multiple edges. If $H$ has a 3-cycle then it has no distinct multiple edges. In these cases we only have $(2,n)$-torus knots and links and twist knots. \hfill$\Box$


\subsection{Links that are s-majors of the $(2,5)$-torus knot}

\begin{figure}[http]
\begin{center}
\includegraphics[width=68mm]{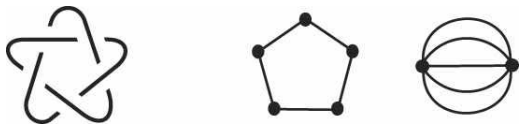}
\caption{The $(2,5)$-torus knot and its graphs}
\label{fig0207}
\end{center}
\end{figure}

\begin{figure}[http]
\begin{center}
\includegraphics[width=50mm]{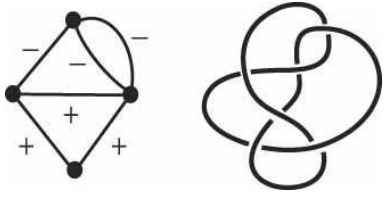}
\caption{A graph corresponding to a $(2,5)$-torus knot}
\label{fig02071}
\end{center}
\end{figure}

A diagram of the $(2,5)$-torus knot and plane graphs $C_5$ and $\theta_5$ corresponding to its underlying projection is illustrated in Fig. \ref{fig0207}.
A graph shown in Fig. \ref{fig02071} below also corresponds to a $(2,5)$-torus knot. It will appear in the following proof. 

Before we describe our next result, we shall confirm some links and their corresponding graphs which are appeared in the statement. 

\begin{figure}[http]
\begin{center}
\includegraphics[width=66mm]{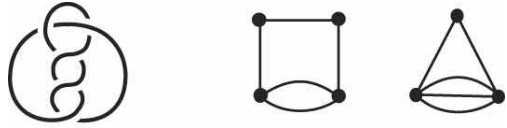}
\caption{$5_2$ and its graphs}
\label{fig020721}
\end{center}
\end{figure}

A graph corresponding to $5_2$ is a graph consists of a 4-cycle with multiple edges. Its dual is a graph consists of a 3-cycle with multiple three edges. See Fig. \ref{fig020721}.

\begin{figure}[http]
\begin{center}
\includegraphics[width=45mm]{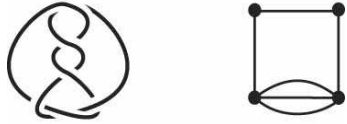}
\caption{$6^2_2$ and its graph}
\label{fig020722}
\end{center}
\end{figure}

A graph corresponding to $6^2_2$ is a graph consists of a 4-cycle with multiple three edges. See Fig. \ref{fig020722}.

\begin{figure}[http]
\begin{center}
\includegraphics[width=67mm]{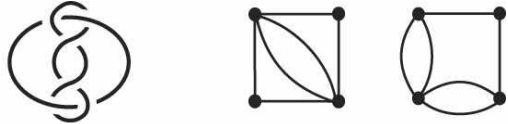}
\caption{$6^2_3$ and its graphs}
\label{fig020723}
\end{center}
\end{figure}

\begin{figure}[http]
\begin{center}
\includegraphics[width=34mm]{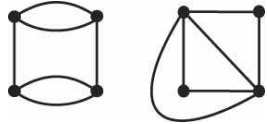}
\caption{Graphs corresponding to $6^2_3$}
\label{fig020723-2}
\end{center}
\end{figure}

A graph corresponding to $6^2_3$ is a graph consists of a 4-cycle with one diagonal multiple edges. Its dual consists of a 4-cycle with two adjacent multiple edges. See Fig. \ref{fig020723}. Graphs illustrated in Fig. \ref{fig020723-2} below also represents $6^2_3$, which is appeared in the following proof. 

\begin{figure}[http]
\begin{center}
\includegraphics[width=70mm]{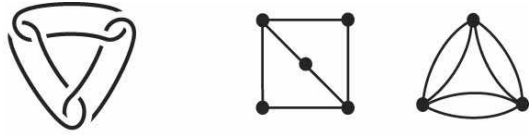}
\caption{$6^3_1$ and its graphs}
\label{fig020724}
\end{center}
\end{figure}

A graph corresponding to $6^3_1$ is the complete bipartite graph $K_{2,3}$. Its dual is a 3-cycle with three multiple edges. See Fig. \ref{fig020724}.

\begin{figure}[http]
\begin{center}
\includegraphics[width=50mm]{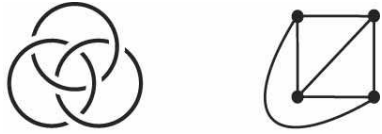}
\caption{$6^3_2$ and its graph}
\label{fig020725}
\end{center}
\end{figure}

A graph corresponding to the Borromean rings $6^3_2$ is the complete graph $K_4$. See Fig. \ref{fig020725}.

\begin{figure}[http]
\begin{center}
\includegraphics[width=68mm]{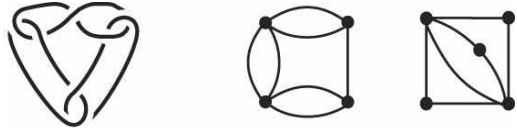}
\caption{$7^3_1$ and its graphs}
\label{fig020726}
\end{center}
\end{figure}

A graph corresponding to $7^3_1$ consists of a 4-cycle with three multiple edges. Its dual is the graph obtained from the complete bipartite graph $K_{2,4}$ by an edge contraction. See Fig. \ref{fig020726}.

\begin{figure}[http]
\begin{center}
\includegraphics[width=70mm]{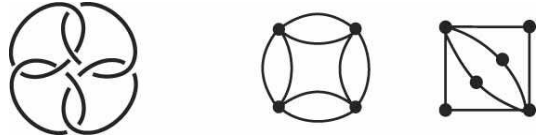}
\caption{$8^4_1$ and its graphs}
\label{fig020727}
\end{center}
\end{figure}

A graph corresponding to $8^4_1$ consists of a 4-cycle with four multiple edges. Its dual is the complete bipartite graph $K_{2,4}$. See Fig. \ref{fig020727}.

{\bf Theorem 2.9.}\quad
{\it A link is an s-major of the $(2,5)$-torus knot if and only if it is non-splittable and it has a prime factor whose crossing number is greater than four that is none of $5_2$, $5^2_1$, $6^2_2$, $6^2_3$, $6^3_1$, $6^3_2$, $7^3_1$, and $8^4_1$.}

{\it Proof.}\quad
It is easily seen that a connected sum of some of knots and links of crossing number less than five, $5_2$, $5^2_1$, $6^2_2$, $6^2_3$, $6^3_1$, $6^3_2$, $7^3_1$, and $8^4_1$ is not an s-major of the $(2,5)$-torus knot. Therefore the \lq only if' part follows.
Let $L$ be a non-splittable link that has a prime factor with crossing number greater than four which is none of $5_2$, $5^2_1$, $6^2_2$, $6^2_3$, $6^3_1$, $6^3_2$, $7^3_1$, and $8^4_1$. Let $\hat L$ be a reduced projection of $L$ and let ${\cal C}$ be a coloring of ${\mathbb S}^2-\hat L$. By the assumption there is a block $H$ of $G(\hat L,{\cal C})$ that corresponds to not necessarily one but some prime factors of $L$ such that at least one of them has crossing number greater than four and it is none of $5_2$, $5^2_1$, $6^2_2$, $6^2_3$, $6^3_1$, $6^3_2$, $7^3_1$, and $8^4_1$.

If $H$ has a cycle of length greater than or equal to five, or if $H$ has two vertices connected by five or more internally disjoint paths, then we obtain graphs corresponding to a $(2,5)$-torus knot. 

If $H$ has a 4-cycle with multiple edges and a diagonal, then we obtain the graph corresponding to a $(2,5)$-torus knot.

Thus we may assume that all cycles of $H$ have length at most four, that every pair of vertices has at most four internally disjoint paths between them, and that if $H$ has a 4-cycle then it does not have multiple edges and diagonals simultaneously.

\begin{figure}[http]
\begin{center}
\includegraphics[width=47mm]{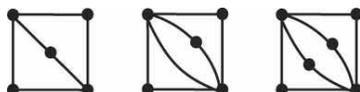}
\caption{4-cycle with a diagonal path of length 2}
\label{4cycle2diagonal}
\end{center}
\end{figure}

If $H$ has a 4-cycle with a diagonal path of length two, such graphs are only three types shown in Fig. \ref{4cycle2diagonal}.

\begin{figure}[http]
\begin{center}
\includegraphics[width=68mm]{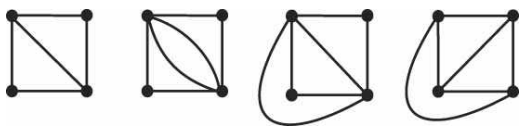}
\caption{4-cycle with a diagonal}
\label{4cycle1diagonal}
\end{center}
\end{figure}

If $H$ has a 4-cycle with a diagonal, such graphs are only four types shown in Fig. \ref{4cycle1diagonal}.

\begin{figure}[http]
\begin{center}
\includegraphics[width=104mm]{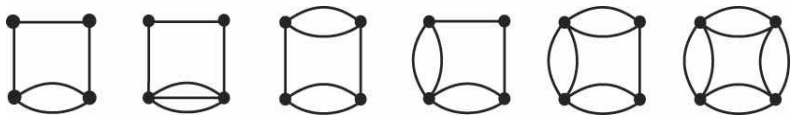}
\caption{4-cycle with multiple edges}
\label{4cycle-multiple}
\end{center}
\end{figure}

If $H$ has a 4-cycle with multiple edges, such graphs are only six types shown in Fig. \ref{4cycle-multiple}. 

\begin{figure}[http]
\begin{center}
\includegraphics[width=50mm]{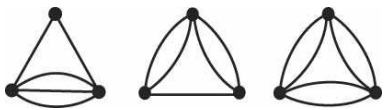}
\caption{3-cycle with multiple edges}
\label{3cycle}
\end{center}
\end{figure}

If $H$ has a 3-cycle with multiple edges, such graphs are only three types shown in Fig. \ref{3cycle}. This completes the proof. \hfill$\Box$

\subsection{Links that are s-majors of $5_2$}

{\bf Theorem 2.10.}\quad
{\it A link is an s-major of $5_2$ if and only if it is non-splittable and it has a prime factor with crossing number greater than four that is none of $(2,n)$-torus knots and links, $5^2_1$, $6^3_1$, and $6^3_2$.}

{\it Proof.}\quad
It is easily seen that a connected sum of some of knots and links of crossing number less than five, $(2,n)$-torus knots and links, $5^2_1$, $6^3_1$, and $6^3_2$ is not an s-major of $5_2$. Therefore the \lq only if' part follows.
Let $L$ be a non-splittable link that has a prime factor with crossing number greater than four which is none of $(2,n)$-torus knots and links, $5^2_1$, $6^3_1$, and $6^3_2$. Let $\hat L$ be a reduced projection of $L$ and let ${\cal C}$ be a coloring of ${\mathbb S}^2-\hat L$. By the assumption there is a block $H$ of $G(\hat L,{\cal C})$ that corresponds to not necessarily one but some prime factors of $L$ such that at least one of them has crossing number greater than four and it is none of $(2,n)$-torus knots and links, $5^2_1$, $6^3_1$, and $6^3_2$.

If $H$ has two vertices which are connected by three internally disjoint paths one of which is of length at least three, or if $H$ has two vertices which are connected by a path of length two and by another at least three internally disjoint paths, then we can get a graph of $5_2$.  

\begin{figure}[http]
\begin{center}
\includegraphics[width=124mm]{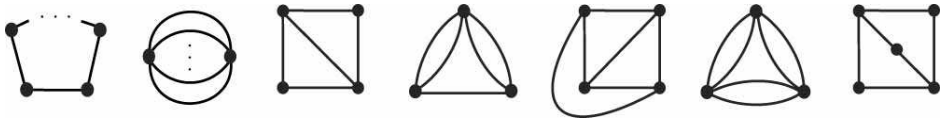}
\caption{Graphs that satisfies the conditions in the proof of Theorem 2.10}
\label{fig02081}
\end{center}
\end{figure}

Thus, we may assume that every pair of vertices has one path of length at least three and one another path, that every pair of vertices has one path of length two and at least two other paths, or that every pair of vertices has paths of length one. In this case, the complete list of such graphs is given in Fig. \ref{fig02081} by the similar way as the previous section, and we get graphs of the $(2,n)$-torus link, $5^2_1$, $6^3_1$, and $6^3_2$. \hfill$\Box$

\section{Final Remark}

In our proofs, the smoothing operation for link diagrams, in other words, the edge-contraction or deletion for the corresponding plane graphs takes an important role. We note a question on the converse of Proposition 1.2.

{\bf Question 3.1.}\quad
{\it Let $L_1$ and $L_2$  be links of the same number of components. Is $L_1\succeq L_2$ implies $L_1 \geq L_2$?}

\end{document}